\documentclass[12pt,reqno]{amsart}
\usepackage{amsfonts,amsmath,amssymb}
\usepackage{graphicx}
\usepackage{psfrag}

\allowdisplaybreaks
\advance\textheight by \topskip

\newcommand{\gauss}[2]{\genfrac{[}{]}{0pt}{}{#1}{#2}_q}
\newcommand{\gausss}[2]{\genfrac{[}{]}{0pt}{}{#1}{#2}_{1/q}}

\numberwithin{equation}{section}

\def\ii{\mathbf{i}}

\begin{document}

\title[Touchard's continued fraction]
{On Touchard's continued fraction and extensions: combinatorics-free, self-contained proofs}
 
\author[H. Prodinger]{Helmut Prodinger}
\address[H. Prodinger]{Department of Mathematics\\
University of Stellenbosch\\
7602 Stellenbosch\\
South Africa}
\email{hproding@sun.ac.za}
\thanks{The author is supported by an incentive grant from the South African National Research Foundation}

\dedicatory{Dedicated to Philippe Flajolet}
\keywords{}
\subjclass{}
\date{\today}
\begin{abstract}
We give a direct and simple proof of Touchard's continued fraction, provide an
extension of it, and transform it into similar expansions related to Motzkin and
Schr\"oder numbers. Another proof is then given that uses only induction.
We use this machinery on two examples that appear in recent papers of
Josuat-Verg\`es; with an additional parameter, these two can be treated simultaneously.
\end{abstract}

\maketitle

\section{Introduction}

Touchard~\cite{Touchard52} studied a certain function $F$ that he developed as a continued fraction (equation (20), loc. cit.):
\begin{equation*}
\frac{1}{1-z(1-q)F(q;z(1-q))}=
\cfrac{1}{1-\cfrac{(1-q)z}{1-\cfrac{(1-q^2)z}{\dots}}}.
\end{equation*}
He (equation (32), loc. cit.) also found the identity
\begin{equation*}
\frac{1}{1-z(1-q)F(q;z(1-q))}=
\frac1{1-zC(z)}\sum_{k\ge0}q^{\binom{k+1}2}\bigl(1-C(z)\bigr)^k,
\end{equation*}
where $C(z)$ is the generating function of the Catalan numbers:
\begin{equation*}
C(z)=\frac{1-\sqrt{1-4z}}{2z}=\sum_{n\ge0}\frac1{n+1}\binom{2n}nz^n.
\end{equation*}
Thanks to the substitution $z=\frac v{(1+v)^2}$, which implies $C(z)=1+v$, 
\emph{Touchard's continued fraction} takes the following form (a ``T-fraction''):
 \begin{equation}\label{TCF}
\cfrac{1}{1+v-\cfrac{(1-q)v}{1+v-\cfrac{(1-q^2)v}{\dots}}}=
\sum_{k\ge0}(-1)^kq^{\binom{k+1}2}v^{k}.
\end{equation}
Touchard's combinatorial approach (counting chord diagrams) is ingenious but at the
same time quite involved. The aim of the present note is to provide a simple and direct proof of (\ref{TCF}).

\medskip

This project was started in June 2009 and was planned to be a cooperation between the present author and Philippe Flajolet. Unfortunately, this never materialised, and it was continued in 2011
with more continued fractions, taken from papers by Josuat-Verg\`es~\cite{josuat1, josuat2}.

There is apparently a rich combinatorial world behind these continued fractions.
However, here, the emphasis is to leave the combinatorics completely out and do
``purely manipulative'' proofs. 

Perhaps the simplest proof is by induction, if one ``knows'' already the result for the
continued fraction with the first $i$ lines deleted. This point of view is elaborated on 
two examples. 

We did not go for an exhaustive list of continued fractions of the Touchard type,
in order to keep this paper short and crisp. However, we are confident that the
approach(es) taken here would also work in other similar instances.

This paper is not overlap-free from other papers, in particular \cite{josuat1, josuat2},
but we are confident that it contains enough original material.

\section{The proof} 
We introduce a parameter $t$ and consider
\begin{equation*}
F(t):=\cfrac{1}{1+v-\cfrac{(1-qt)v}{1+v-\cfrac{(1-q^2t)v}{\dots}}}=
\frac{1}{1+v-(1-qt)vF(qt)}.
\end{equation*}
Setting $F(t)=A(t)/B(t)$, this leads to $A(t)=B(qt)$ and
\begin{equation*}
B(t)=(1+v)B(qt)-(1-qt)vB(q^2t).
\end{equation*}
We make the \emph{ansatz} 
\begin{equation*}
B(t)=\sum_{n\ge0}a_nt^n
\end{equation*}
and get upon comparing coefficients
\begin{equation*}
a_n=(1+v)q^na_n-vq^{2n}a_n+vq^{2n-1}a_{n-1}=
\frac{vq^{2n-1}a_{n-1}}{(1-q^n)(1-vq^n)}.
\end{equation*}
This can be iterated, and since $a_0=1$, 
\begin{equation*}
a_n=\frac{v^nq^{n^2}}{(q;q)_n(vq;q)_n},
\end{equation*}
where we employed the notation $(x;q)_n=(1-x)(1-xq)\dots(1-xq^{n-1})$.
Thus 
\begin{equation*}
B(t)=\sum_{n\ge0}\frac{v^nq^{n^2}}{(q;q)_n(vq;q)_n}t^n,
\end{equation*}
and \emph{Touchard's continued fraction} is given by $F(1)=B(q)/B(1)$.

Our goal will be achieved once we are able to establish the identity
\begin{equation*}
B(q)=B(1)\sum_{k\ge0}(-1)^kq^{\binom{k+1}2}v^{k}.
\end{equation*}
This is the same as
\begin{equation*}
\sum_{n\ge0}\frac{v^nq^{n^2}}{(q;q)_n(v;q)_n}=
\sum_{n\ge0}\frac{v^nq^{n^2-n}}{(q;q)_n(v;q)_n}\cdot\sum_{n\ge0}(-1)^nq^{\binom{n}2}v^{n},
\end{equation*}
where we replaced $v$ by $v/q$. One of the sums can be evaluated thanks to a formula due to Cauchy~\cite[(10.9.2)]{AnAsRo99}:
\begin{equation*}
\sum_{n\ge0}\frac{v^nq^{n^2-n}}{(q;q)_n(v;q)_n}
=\frac1{(v;q)_\infty}=\sum_{n\ge0}\frac{v^n}{(q;q)_n}.
\end{equation*}
Comparing coefficients of $v^n$, we are left to prove that%
\footnote{J.~Cigler, who has seen an early draft of this paper,
has shown me a ``human'' proof.}
\begin{equation*}
\sum_{k=1}^n\frac{q^{k^2}}{(q;q)_k}\gauss{n-1}{k-1}=
\sum_{k=0}^n\frac{(-1)^kq^{\binom k2}}{(q;q)_{n-k}}.
\end{equation*}
Here, we use the notation
\begin{equation*}
\gauss{n}{k}=\frac{(q;q)_{n}}{(q;q)_{k}(q;q)_{n-k}}.
\end{equation*}
To establish such an identity is nowadays routine, and the tool
for it is the $q$-version of Zeilberger's algorithm~\cite{PeWiZe96}.
It  produces the recursion 
\begin{equation*}
(q^{n+1}-q)T(n)+(q^{2n+3}-q^{n+2}+q-q^{n+1}+1)T(n+1)+(q^{n+2}-1)T(n+2)=0,
\end{equation*}
for both sides, and together with a few easily checked initial conditions, we are done.

\section{An extension}

Our proof gives more, by considering $F(q^i)$:
\begin{equation*}
F(q^i)=\frac{B(q^{i+1})}{B(q^i)}=\cfrac{1}{1+v-\cfrac{(1-q^i)v}{1+v-\cfrac{(1-q^{i+1})v}{\dots}}}.
\end{equation*}
This quotient by itself it not nice, but we have
\begin{equation*}
\frac{B(q^{i+1})}{B(1)}=\sum_{n\ge0}(-1)^n\gauss{n+i}{i}q^{\binom{n+1}2}v^n.
\end{equation*}
The identity behind this is
\begin{equation*}
\sum_{n\ge0}\frac{q^{n^2+in}v^n}{(q)_n(v)_n}=
\sum_{n\ge0}\frac{v^n}{(q)_n}\cdot
\sum_{n\ge0}(-1)^nv^nq^{\binom n2}\gauss{n+i}{i},
\end{equation*}
which, in terms of coefficients, reads as
\begin{equation*}
\sum_{k=1}^n\frac{q^{k^2+ik}}{(q;q)_k}\gauss{n-1}{k-1}=
\sum_{k=0}^n\frac{(-1)^kq^{\binom k2}}{(q;q)_{n-k}}\gauss{k+i}{i}.
\end{equation*}
Again, Zeilberger's algorithm proves this readily. Therefore
\begin{align*}
F(q^i)&=\frac{B(q^{i+1})}{B(1)}\bigg/\frac{B(q^{i})}{B(1)}\\*&
=\sum_{n\ge0}(-1)^n\gauss{n+i}{i}q^{\binom{n+1}2}v^n\bigg/
\sum_{n\ge0}(-1)^n\gauss{n+i-1}{i-1}q^{\binom{n+1}2}v^n.
\end{align*}

\section{Another proof, starting from the righthand side}

We set 
\begin{equation*}
s_i:=\sum_{n\ge0}(-1)^n\gauss{n+i}{i}q^{\binom{n+1}2}v^n\quad i\ge0,
\end{equation*}
$s_{-1}:=1$
and expand the quantity $\dfrac{s_0}{s_{-1}}$ as a continued fraction of the form
\begin{equation*}
\cfrac{1}{1+v-\cfrac{a_1v}{1+v-\cfrac{a_2v}{\dots}}}.
\end{equation*}
It is not hard to show that the sequence $a_1, a_2,\dots$ is unique, as one can compare
coefficients in series expansions of both sides, and compute them one by one. 

We claim that $a_i=1-q^i$ and prove the formula
\begin{equation*}
\cfrac{1}{1+v-\cfrac{a_{i+1}v}{1+v-\cfrac{a_{i+2}v}{\dots}}}=\frac{s_{i}}{s_{i-1}}
\end{equation*}
by induction. Since $s_0=1-qv+\cdots$, the value $a_1$ is established. And now
\begin{equation*}
\frac{1}{1+v-a_{i+1}v\dfrac{s_{i+1}}{s_{i}}}=\frac{s_{i}}{s_{i-1}}
\end{equation*}
leads to
\begin{equation}\label{recu1}
(1+v)s_i-a_{i+1}vs_{i+1}=s_{i-1}.
\end{equation}
(All the steps are reversible, and it might be clearer to start the arguments from this
recursion.)
Comparing coefficients,
\begin{align*}
(-1)^n\gauss{n+i}{i}q^{\binom{n+1}2}&+(-1)^{n-1}\gauss{n+i-1}{i}q^{\binom{n}2}
-(1-q^{i+1})(-1)^{n-1}\gauss{n+i}{i+1}q^{\binom{n}2}\\&=
(-1)^n\gauss{n+i-1}{i-1}q^{\binom{n+1}2},
\end{align*}
or equivalently
\begin{align*}
\gauss{n+i}{i}q^{n}&-\gauss{n+i-1}{i}
+(1-q^{i+1})\gauss{n+i}{i+1}=
\gauss{n+i-1}{i-1}q^{n},
\end{align*}
or
\begin{align*}
\gauss{n+i}{i}q^{n}&-\gauss{n+i-1}{i}
+(1-q^{n+i})\gauss{n+i-1}{i}=
\gauss{n+i-1}{i-1}q^{n},
\end{align*}
or
\begin{align*}
\gauss{n+i}{i}-q^{i}\gauss{n+i-1}{i}=
\gauss{n+i-1}{i-1},
\end{align*}
which is the basic recursion; thus the proof is finished.

\subsection*{Remark.} J.~Cigler has kindly pointed out that
\begin{equation*}
s_i=\sum_{k\ge0}\gausss{-i-1}{k}v^k,
\end{equation*}
which are Rogers-Szeg\H{o} polynomial with negative indices.

Furthermore, for $v=-1$ \eqref{recu1} simplifies, and the recursion can be solved by iteration,
leading to a product representation.

\section{Motzkin numbers}

Cigler~\cite{private_communication} recently\footnote{In 2009, when a first draft of this paper was sketched.} found	the identity
\begin{equation*}
\cfrac{1}{1-z-\cfrac{(1-q)z^2}{1-z-\cfrac{(1-q^2)z^2}{\dots}}}=
\sum_{k\ge0}(-1)^kq^{\binom{k+1}2}z^{2k}M(z)^{2k+1},
\end{equation*}
where
\begin{equation*}
M(z)=\frac{1-z-\sqrt{1-2z-3z^2}}{2z^2}
\end{equation*}
is the generating function of the Motzkin numbers.

Using the substitution $z=\frac{v}{1+v+v^2}$, we see that, after division of both
sides by $1+v+v^2$, only even powers of
$v$ appear, and with the further substitution $v^2=x$, the formula transforms
into \emph{Touchard's continued fraction}.

\section{Schr\"oder numbers}

Cigler also found the formula
\begin{equation*}
\cfrac{1}{1-z-\cfrac{(1-q)z}{1-z-\cfrac{(1-q^2)z}{\dots}}}=
\sum_{k\ge0}(-1)^kq^{\binom{k+1}2}z^{k}S(z)^{2k+1},
\end{equation*}
where
\begin{equation*}
S(z)=\frac{1-z-\sqrt{1-6z+z^2}}{2z}
\end{equation*}
is the generating function of the Schr\"oder numbers.

We can also reduce this to an instance of \emph{Touchard's continued fraction}
as follows: The substitution $y=z/(1-z)^2$ transforms the formula into
\begin{equation*}
\cfrac{1}{1-\cfrac{(1-q)y}{1-\cfrac{(1-q^2)y}{\dots}}}=
(1-z)\sum_{k\ge0}(-1)^kq^{\binom{k+1}2}z^{k}S(z)^{2k+1}.
\end{equation*}
But recall that
\begin{equation*}
\cfrac{1}{1-\cfrac{(1-q)y}{1-\cfrac{(1-q^2)y}{\dots}}}=
\frac1{1-yC(y)}\sum_{k\ge0}q^{\binom{k+1}2}\bigl(1-C(y)\bigr)^k,
\end{equation*}
so that the claim  follows from the equalities
\begin{equation*}
(1-z)S(z)=\frac1{1-yC(y)}\quad\text{and}\quad
zS^2(z)=C(y)-1.
\end{equation*}
But both are a routine verification that is best done with a computer algebra system.

\section{A formula derived by Riordan}
Writing
\begin{equation*}
\frac{1}{1-zF(q;z)}=\sum_{n\ge0}T_n(q)z^n,
\end{equation*}
we find, upon comparing coefficients, that
\begin{align*}
T_n(q)(1-q)^n&=[z^n]
\frac1{1-zC(z)}\sum_{k\ge0}q^{\binom{k+1}2}\bigl(1-C(z)\bigr)^k\\
&=[v^n](1-v)(1+v)^{2n}\sum_{k\ge0}q^{\binom{k+1}2}(-v)^k\\
&=\sum_{k=0}^n(-1)^kq^{\binom{k+1}2}\biggl[
\binom{2n}{n-k}-\binom{2n}{n-k-1}\biggr],
\end{align*}
which is a formula derived by Riordan~\cite{Riordan75}.

\section{Another continued fraction}

We consider now an example from~\cite{josuat1, josuat2}:
\begin{equation*}
\cfrac{1}{1+v-\cfrac{(1-q)^2v}{1+v-\cfrac{(1-q^2)^2v}{\dots}}}.
\end{equation*}
We set
\begin{equation*}
F(t)=\cfrac{1}{1+v-\cfrac{(1-qt)^2v}{1+v-\cfrac{(1-q^2t)^2v}{\dots}}}
=\frac1{1+v-(1-qt)^2vF(qt)}.
\end{equation*}
This continued fractions appears starting from 
\begin{equation*}
\cfrac{1}{1-\cfrac{(1-q)^2z}{1-\cfrac{(1-q^2)^2z}{\dots}}},
\end{equation*}
after the usual substitution. Set 
\begin{equation*}
F(t)=\frac{A(t)}{B(t)}.
\end{equation*}
We get $A(t)=B(qt)$ and 
\begin{equation*}
B(t)=(1+v)B(qt)-(1-qt)^2vB(q^2t).
\end{equation*}
Reading off coefficients of $t^n$ leads to a recursion of second order, which is not nice.
But if we define
\begin{equation*}
\beta(t):=\frac{B(t)}{(t;q)_\infty},
\end{equation*}
then
\begin{equation*}
(1-t)\beta(t)=(1+v)\beta(qt)-(1-qt)v\beta(q^2t),
\end{equation*}
and
\begin{equation*}
\beta_n-\beta_{n-1}=(1+v)q^n\beta_n-vq^{2n}\beta_n+q^{2n-1}v\beta_{n-1},
\end{equation*}
or
\begin{equation*}
\beta_n=\beta_{n-1}\frac{1+vq^{2n-1}}{(1-q^n)(1-vq^n)}=\frac{(-vq;q^2)_n}{(q;q)_n(vq;q)_n},
\end{equation*}
since $\beta_0=1$. Hence
\begin{equation*}
\beta(t)=\sum_{n\ge0}\frac{(-vq;q^2)_n}{(q;q)_n(vq;q)_n}t^n
\end{equation*}
and
\begin{equation*}
F(t)=\frac{\beta(qt)}{(1-t)\beta(t)}.
\end{equation*}
We need $F(1)$. The interpretation of $(1-t)\beta(t)$ is as a limit. We use the Heine transform,
as was done in ~\cite{josuat1, josuat2} as well:
\begin{equation*}
\beta(t)=\sum_{n\ge0}\frac{(\ii \sqrt{vq};q)_n(-\ii \sqrt{vq};q)_n}{(q;q)_n(vq;q)_n}t^n,
\end{equation*}
and therefore
\begin{equation*}
\beta(t)=\frac{(\ii \sqrt{vq}t;q)_\infty(-\ii \sqrt{vq};q)_\infty}{(vq;q)_\infty(t;q)_\infty}
\sum_{n\ge0}\frac{(\ii \sqrt{vq};q)_n(t;q)_n}{(q;q)_n(\ii \sqrt{vq}t;q)_n}(-\ii \sqrt{vq}\,)^n
\end{equation*}
and
\begin{equation*}
(1-t)\beta(t)=\frac{(\ii \sqrt{vq}t;q)_\infty(-\ii \sqrt{vq};q)_\infty}{(vq;q)_\infty(tq;q)_\infty}
\sum_{n\ge0}\frac{(\ii \sqrt{vq};q)_n(t;q)_n}{(q;q)_n(\ii \sqrt{vq}t;q)_n}(-\ii \sqrt{vq}\,)^n.
\end{equation*}
Now we can plug in $t=1$:
\begin{equation*}
\frac{(-vq;q^2)_\infty}{(vq;q)_\infty(q;q)_\infty}.
\end{equation*}
Furthermore
\begin{equation*}
\beta(q)=\frac{(\ii \sqrt{vq}q;q)_\infty(-\ii \sqrt{vq};q)_\infty}{(vq;q)_\infty(q;q)_\infty}
\sum_{n\ge0}\frac{(\ii \sqrt{vq};q)_n(q;q)_n}{(q;q)_n(\ii \sqrt{vq}q;q)_n}(-\ii \sqrt{vq}\,)^n,
\end{equation*}
or
\begin{equation*}
\beta(q)=\frac{(-vq;q^2)_\infty}{(vq;q)_\infty(q;q)_\infty}
\sum_{n\ge0}\frac{1}{1-\ii \sqrt{vq}q^n}(-\ii \sqrt{vq}\,)^n.
\end{equation*}
Hence
\begin{align*}
F(1)&=\sum_{n\ge0}\frac{1}{1-\ii \sqrt{vq}q^n}(-\ii \sqrt{vq}\,)^n\\
&=\sum_{n,m\ge0}(\ii \sqrt{vq}q^n)^m(-\ii \sqrt{vq}\,)^n\\
&=\sum_{n,m\ge0}\ii ^{m+n}(\sqrt{vq}\,)^{m+n}q^{nm}(-1)^n\\
&=\sum_{0\le n\le 2N}(-1)^{N}(vq)^{N}q^{n(2N-n)}(-1)^n.
\end{align*}
So
\begin{align*}
[v^N]F(1)&=\sum_{0\le n\le 2N}(-1)^{N+n}q^{N+2nN-n^2}\\
&=\sum_{-N\le n\le N}(-1)^{n}q^{N+2(n+N)N-(n+N)^2}\\
&=q^{N(N+1)}\sum_{-N\le n\le N}(-1)^{n}q^{-n^2}.
\end{align*}
To find the coefficients (in the variable $z$) of the original continued fraction,
we compute
\begin{align*}
[z^N](1+v)F(1)&=[v^N](1+v)^{2N}(1-v)F(1)\\
&=\sum_{k=0}^{N}[v^{N-k}](1+v)^{2N}[v^k](1-v)F(1)\\
&=\sum_{k=0}^{N}\biggl[\binom{2N}{N-k}-\binom{2N}{N-k-1}\biggr]
q^{k(k+1)}\sum_{-k\le j\le k}(-1)^{j}q^{-j^2}.
\end{align*}
This translation from the coefficients in the variable $v$ to $z$ works always
in exactly the same way, so we will not state the $z$-formula for the further
examples. 

\medskip

As a bonus, we can also consider
\begin{equation*}
\cfrac{1}{1+v-\cfrac{(1-q^i)^2v}{1+v-\cfrac{(1-q^{i+1})^2v}{\dots}}},
\end{equation*}
which can be computed via
\begin{equation*}
\frac{B(q^{i+1})}{B(q^{i})}=\frac{B(q^{i+1})}{B(1)}\bigg/\frac{B(q^{i})}{B(1)}.
\end{equation*}
Now
\begin{align*}
\frac{B(q^{i})}{B(1)}&=\frac{(q^i;q)_\infty\beta(q^i)}{\frac{(-vq;q^2)_\infty}{(vq;q)_\infty}}\\
&=\frac{1}{(\ii \sqrt{vq};q)_i}
\sum_{n\ge0}\frac{(\ii \sqrt{vq};q)_n(q^i;q)_n}{(q;q)_n(\ii \sqrt{vq}q^i;q)_n}(-\ii \sqrt{vq}\,)^n\\
&=\sum_{n\ge0}\frac{(\ii \sqrt{vq};q)_n(q^i;q)_n}{(q;q)_n(\ii \sqrt{vq};q)_{n+i}}(-\ii \sqrt{vq}\,)^n\\
&=\frac1{(q;q)_{i-1}}\sum_{n\ge0}\frac{(q^{n+1};q)_{i-1}}{(\ii \sqrt{vq}q^n;q)_{i}}(-\ii \sqrt{vq}\,)^n\\
&=\frac1{(q;q)_{i-1}}\sum_{n,m\ge0}\gauss{i+m-1}{m}(q^{n+1};q)_{i-1}
(\ii \sqrt{vq}q^n)^m(-\ii \sqrt{vq}\,)^n\\
&=\frac1{(q;q)_{i-1}}\sum_{N\ge0}v^N\sum_{0\le n\le 2N}\gauss{i+2N-n-1}{i-1}(q^{n+1};q)_{i-1}
q^{N+2nN-n^2}(-1)^{N+n}\\
&=\frac1{(q;q)_{i-1}}\sum_{N\ge0}v^Nq^{N(N+1)}\sum_{-N\le n\le N}\gauss{i+N-n-1}{i-1}(q^{n+N+1};q)_{i-1}
q^{-n^2}(-1)^{n}\\
&=\sum_{N\ge0}v^Nq^{N(N+1)}\sum_{-N\le n\le N}\gauss{i+N-n-1}{i-1}\gauss{i+N+n-1}{i-1}q^{-n^2}(-1)^{n}\\
\end{align*}
and the continued fraction has been evaluated as a quotient of two series of this type.

\section{Independent proof by induction}
As before we can turn this extra information into an extremely elementary proof by induction:
We set
\begin{equation*}
s_i:=\sum_{N\ge0}v^Nq^{N(N+1)}\sum_{-N\le n\le N}\gauss{i+N-n}{i}\gauss{i+N+n}{i}q^{-n^2}(-1)^{n},
\end{equation*}
$s_{-1}:=1$,
and expand the quantity $\dfrac{s_0}{s_{-1}}$ as a continued fraction of the form
\begin{equation*}
\cfrac{1}{1+v-\cfrac{a_1v}{1+v-\cfrac{a_2v}{\dots}}}.
\end{equation*}
We claim that $a_i=(1-q^i)^2$ and prove the formula
\begin{equation*}
\cfrac{1}{1+v-\cfrac{a_{i+1}v}{1+v-\cfrac{a_{i+2}v}{\dots}}}=\frac{s_{i}}{s_{i-1}}
\end{equation*}
by induction. Since $1/s_0=1-q( q-2 ) v+\cdots$, the value $a_1$ is established. And now
\begin{equation*}
\frac{1}{1+v-a_{i+1}v\dfrac{s_{i+1}}{s_{i}}}=\frac{s_{i}}{s_{i-1}}
\end{equation*}
leads to
\begin{equation}\label{recu2}
(1+v)s_i-a_{i+1}vs_{i+1}=s_{i-1}.
\end{equation}
Set 
\begin{equation*}
s_{i,N}:=q^{N(N+1)}\sum_{-N\le n\le N}\gauss{i+N-n}{i}\gauss{i+N+n}{i}q^{-n^2}(-1)^{n}.
\end{equation*}
Then we must show the above recursion, which, in terms of the coefficients, reads as
\begin{align*}
s_{i,N}+s_{i,N-1}-(1-q^{i+1})^2s_{i+1,N-1}-s_{i-1,N}=0.
\end{align*}
Manuel Kauers (Risc, Linz) kindly provided an automatic
proof for this, using the package \textsf{HolonomicFunctions} by
Christoph Koutschan~\cite{Koutschan09}.

\subsection*{Remark.}J.~Cigler has pointed out that \eqref{recu2} becomes nicer when written
with negative indices: Set $S_i:=s_{-i-1}$, then
\begin{equation*}
S_i=(1+v)S_{i-1}-(1-q^{1-i})^2vS_{i-2},\qquad S_0=1,\quad S_1=1+v
\end{equation*}
and
\begin{equation*}
S_i=\sum_{N=0}^iv^Nq^{2N^2-2Ni}\sum_{-N\le n\le N}(-1)^n\gauss{i}{N-n}\gauss{i}{N+n}.
\end{equation*}
For $q=1$, this is just the binomial theorem.

\section{Another continued fraction of the Schr\"oder type}

In this section we investigate the following continued fraction: (cf.~\cite[(18)]{josuat1})
\begin{equation*}
\cfrac{1}{1+v-\cfrac{(1-q)(1-q^2)v}{1+v-\cfrac{(1-q^2)(1-q^3)v}{\dots}}}.
\end{equation*}
Set
\begin{equation*}
F(t)=\cfrac{1}{1+v-\cfrac{(1-qt)(1-q^2t)v}{1+v-\cfrac{(1-q^2t)(1-q^3t)v}{\dots}}}
=\frac1{1+v-(1-qt)(1-q^2t)vF(qt)}.
\end{equation*}
With the usual $F(t)=A(t)/B(t)$, we find $A(t)=B(qt)$ and
\begin{equation*}
B(t)=(1+v)B(qt)-(1-qt)(1-q^2t)vB(q^2t).
\end{equation*}
We set $\beta(t)=B(t)/(qt;q)_\infty$:
\begin{equation*}
(1-qt)\beta(t)=(1+v)\beta(qt)-(1-qt)v\beta(q^2t).
\end{equation*}
Hence
\begin{equation*}
\beta_n-q\beta_{n-1}=(1+v)q^n\beta_n-vq^{2n}\beta_n+vq^{2n-1}\beta_{n-1},
\end{equation*}
or
\begin{equation*}
\beta_n=\frac{q(1+vq^{2n-2})}{(1-q^n)(1-vq^n)}\beta_{n-1}
=\frac{q^n(-v;q^2)_n}{(q;q)_n(vq;q)_n}.
\end{equation*}
So
\begin{align*}
\beta(t)&=\sum_{n\ge0}\frac{(\ii \sqrt v;q)_n(-\ii \sqrt v;q)_n}{(q;q)_n(vq;q)_n}(tq)^n\\
&=\frac{(-\ii \sqrt v;q)_\infty(\ii\sqrt vqt;q)_\infty}{(vq;q)_\infty(qt;q)_\infty}
\sum_{n\ge0}\frac{(\ii\sqrt v q;q)_n(qt;q)_n}{(q;q)_n(\ii\sqrt vqt;q)_n}(-\ii \sqrt v\,)^n.
\end{align*}
Hence
\begin{align*}
\beta(1)
&=\frac{(-v;q^2)_\infty}{(vq;q)_\infty(q;q)_\infty}
\frac1{1+\ii\sqrt v}
\end{align*}
and
\begin{align*}
\beta(q)
&=\frac{(-\ii \sqrt v;q)_\infty(\ii\sqrt vq^2;q)_\infty}{(vq;q)_\infty(q^2;q)_\infty}
\sum_{n\ge0}\frac{(\ii\sqrt v q;q)_n(q^2;q)_n}{(q;q)_n(\ii\sqrt vq^2;q)_n}(-\ii \sqrt v\,)^n\\
&=\frac{(-v;q^2)_\infty}{(vq;q)_\infty(q;q)_\infty}
\sum_{n\ge0}\frac{(1-q^{n+1})}{(1-\ii\sqrt vq^{n+1})}(-\ii \sqrt v\,)^n.
\end{align*}
Therefore
\begin{align*}
F(1)&=\frac{B(q)}{B(1)}=\frac{(q^2;q)_\infty\beta(q)}{(q;q)_\infty\beta(1)}\\
&=\frac{1+\ii\sqrt v}{1-q}\sum_{n\ge0}\frac{(1-q^{n+1})}{(1-\ii\sqrt vq^{n+1})}(-\ii \sqrt v\,)^n\\
&=\frac{1+\ii\sqrt v}{1-q}\sum_{n,m\ge0}(1-q^{n+1})(\ii\sqrt vq^{n+1})^m(-\ii \sqrt v\,)^n\\
&=\frac{1}{1-q}\sum_{N\ge0}\sum_{0\le n\le 2N}(1-q^{n+1})(-1)^{N+n}v^Nq^{(n+1)(2N-n)}\\
&\quad+\frac{1}{1-q}\sum_{N\ge0}\sum_{0\le n\le 2N+1}(1-q^{n+1})(-1)^{N+n+1}v^{N+1}q^{(n+1)(2N+1-n)}\\
&=\frac{1}{1-q}\sum_{N\ge0}v^Nq^{N(N+1)}\sum_{-N\le n\le N}(1-q^{n+N+1})(-1)^{n}
q^{-n(n+1)}\\
&\quad-\frac{1}{1-q}\sum_{N\ge0}v^{N}q^{N^2}\sum_{-N\le n\le N}(1-q^{n+N})(-1)^{n}q^{-n^2}.
\end{align*}
One could rearrange that, but it would not get any better.

\section{Generalisation of the last two continued fractions}

Here we consider
\begin{equation*}
\cfrac{1}{1+v-\cfrac{(1-q)(1-q^d)v}{1+v-\cfrac{(1-q^2)(1-q^{d+1})v}{\dots}}},
\end{equation*}
with $d$ a nonnegative integer. So for $d=1$ and $d=2$ we get the continued fractions
from the previous sections. As usual
\begin{equation*}
F(t)=\cfrac{1}{1+v-\cfrac{(1-qt)(1-q^dt)v}{1+v-\cfrac{(1-q^2t)(1-q^{d+1}t)v}{\dots}}},
\end{equation*}
so
\begin{equation*}
F(t)=\frac1{1+v-(1-qt)(1-q^dt)vF(qt)}=\frac{B(qt)}{B(q)}
\end{equation*}
with
\begin{equation*}
B(t)=(1+v)B(qt)-(1-qt)(1-q^dt)vB(q^2t).
\end{equation*}
Set $\beta(t)=B(t)/(q^{d-1}t;q)_\infty$, then
\begin{equation*}
(1-q^{d-1}t)\beta(t)=(1+v)\beta(qt)-(1-qt)v\beta(q^2t).
\end{equation*}
Comparing coefficients,
\begin{equation*}
\beta_n-q^{d-1}\beta_{n-1}=(1+v)q^n\beta_n-vq^{2n}\beta_n+vq^{2n-1}\beta_{n-1},
\end{equation*}
or
\begin{equation*}
\beta_n=\frac{q^{d-1}(1+vq^{2n-d})}{(1-q^n)(1-vq^n)}\beta_{n-1}=
\frac{q^{(d-1)n}(-vq^{2-d};q^2)_n}{(q;q)_n(vq;q)_n}.
\end{equation*}
So 
\begin{align*}
\beta(t)&=\sum_{n\ge0}\frac{(\ii \sqrt vq^{1-d/2};q)_n(-\ii \sqrt vq^{1-d/2};q)_n}{(q;q)_n(vq;q)_n}(q^{d-1}t)^n\\
&=\frac{(-\ii \sqrt vq^{1-d/2};q)_\infty(\ii \sqrt vq^{d/2}t;q)_\infty}{(vq;q)_\infty(q^{d-1}t;q)_\infty}
\sum_{n\ge0}\frac{(\ii\sqrt vq^{d/2};q)_n(q^{d-1}t;q)_n}{(q;q)_n(\ii \sqrt vq^{d/2}t;q)_n}
(-\ii \sqrt vq^{1-d/2})^n.
\end{align*}
Let us exclude the case $d=1$, so that we don't have to worry about taking a limit.
\begin{align*}
\beta(1)
&=\frac{(-vq^{2-d};q^2)_\infty}{(vq;q)_\infty(q^{d-1};q)_\infty(\ii \sqrt vq^{1-d/2};q)_{d-1}}
\sum_{n\ge0}\frac{(q^{d-1};q)_n}{(q;q)_n}
(-\ii \sqrt vq^{1-d/2})^n.
\end{align*}
\begin{align*}
\beta(q)
&=\frac{(-vq^{2-d};q^2)_\infty}{(vq;q)_\infty(q^{d};q)_\infty(\ii \sqrt vq^{1-d/2};q)_{d-1}}
\sum_{n\ge0}\frac{(\ii\sqrt vq^{d/2};q)_n(q^{d};q)_n}{(q;q)_n
(\ii \sqrt vq^{1+d/2};q)_{n+1}}
(-\ii \sqrt vq^{1-d/2})^n\\
&=\frac{(-vq^{2-d};q^2)_\infty}{(vq;q)_\infty(q^{d};q)_\infty(\ii \sqrt vq^{1-d/2};q)_{d-1}}
\\&\qquad\qquad\qquad\times\sum_{n\ge0}\frac{(1-\ii\sqrt vq^{d/2})(q^{d};q)_n}{(q;q)_n
(1-\ii \sqrt vq^{1+d/2+n})(1-\ii \sqrt vq^{d/2+n})}
(-\ii \sqrt vq^{1-d/2})^n.
\end{align*}
So
\begin{align*}
F(1)&=\frac{\beta(q)}{(1-q^{d-1})\beta(1)}\\
&=\sum_{n\ge0}\frac{(q^{d};q)_n(1-\ii\sqrt vq^{d/2})}{(q;q)_n
(1-\ii \sqrt vq^{1+d/2+n})(1-\ii \sqrt vq^{d/2+n})}
(-\ii \sqrt vq^{1-d/2})^n\\
&=\sum_{n\ge1}\frac{(q^{d};q)_n}{q^n(1-q)(q;q)_{n-1}
(1-\ii \sqrt vq^{d/2+n})}
(-\ii \sqrt vq^{1-d/2})^n\\
&\quad+\sum_{n\ge0}\frac{(q^{d};q)_n(1-q^{n+1})}{q^n(1-q)(q;q)_n
(1-\ii \sqrt vq^{1+d/2+n})}
(-\ii \sqrt vq^{1-d/2})^n\\
&=\sum_{n\ge1,m\ge0}\frac{(q^{d};q)_n}{q^n(1-q)(q;q)_{n-1}
}(\ii \sqrt vq^{d/2+n})^m
(-\ii \sqrt vq^{1-d/2})^n\\
&\quad+\sum_{n,m\ge0}\frac{(q^{d};q)_n(1-q^{n+1})}{q^n(1-q)(q;q)_n
}(\ii \sqrt vq^{1+d/2+n})^m
(-\ii \sqrt vq^{1-d/2})^n\\
&=\sum_{N\ge0}v^N\sum_{1\le n\le 2N}\frac{(q^{d};q)_n(-1)^{N+n}}{(1-q)(q;q)_{n-1}
} q^{n(2N-n)+Nd-nd}\\
&\quad+\sum_{N\ge0}v^N\sum_{0\le n\le 2N}\frac{(q^{d};q)_n(1-q^{n+1})(-1)^{N+n}}{(1-q)(q;q)_n
}q^{(1+n)(2N-n)+Nd-nd}.
\end{align*}

\subsection*{Acknowledment.} I would like to thank J.~Cigler and P.~Flajolet for
valuable suggestions.

\bibliographystyle{plain}

\end{document}